\documentclass[10pt]{article}
\usepackage{latexsym,amsmath,amssymb,graphics,amscd}

\makeatletter
\@addtoreset{figure}{section}
\def\thefigure{\thesection.\@arabic\c@figure}
\def\fps@figure{h,t}
\@addtoreset{table}{bsection}

\def\thetable{\thesection.\@arabic\c@table}
\def\fps@table{h, t}
\@addtoreset{equation}{section}

\makeatother

\begin{document}

\newtheorem{theorem}{Theorem}[section]
\newtheorem{definition}[theorem]{Definition}
\newtheorem{lemma}[theorem]{Lemma}
\newtheorem{remark}[theorem]{Remark}
\newtheorem{proposition}[theorem]{Proposition}
\newtheorem{corollary}[theorem]{Corollary}
\newtheorem{example}[theorem]{Example}
\newtheorem{examples}[theorem]{Examples}

\newcommand{\bfi}{\bfseries\itshape}

\newsavebox{\savepar}
\newenvironment{boxit}{\begin{lrbox}{\savepar}
\begin{minipage}[b]{15.8cm}}{\end{minipage}
\end{lrbox}\fbox{\usebox{\savepar}}}

\makeatletter
\title{{\bf Cotangent bundle reduction}}
\author{Juan-Pablo Ortega$^{1}$ and  Tudor S. Ratiu$^{2}$}
\addtocounter{footnote}{1}
\footnotetext{Centre National de la Recherche Scientifique,
D\'epartement de Math\'ematiques de Besan\c con,
Universit\'e de Franche-Comt\'e.
UFR des Sciences et Techniques.
16, route de Gray.
F-25030 Besan\c con cedex. France. {\texttt
Juan-Pablo.Ortega@math.univ-fcomte.fr}. }
\addtocounter{footnote}{1}
\footnotetext{Section de Math\'ematiques, 
\'Ecole Polytechnique F\'ed\'erale de Lausanne,  CH-1015 Lausanne,
Switzerland. {\texttt Tudor.Ratiu@epfl.ch}.}
\date{}
\makeatother
\maketitle

\begin{abstract}
This encyclopedia article briefly reviews without proofs some of the main results
in cotangent bundle reduction. The article recalls most the necessary prerequisites to
understand the main results.
\end{abstract}

The general symplectic reduction theory presented in~\cite{or} becomes much richer and has many
applications if the symplectic manifold is the cotangent bundle $(T ^\ast
Q , \Omega_Q = - \mathbf{d}\Theta_Q)$ of a manifold $Q$. The canonical
one-form $\Theta_Q$ on $T ^\ast Q $ is given by $\Theta_Q( \alpha_q)
\left(V_{\alpha_q} \right)= \alpha_q \left(T_{ \alpha_q} \pi_Q
\left(V_{\alpha_q}\right)\right)$, for any $q \in Q , \alpha_q \in
T_q ^\ast Q$, and tangent vector $V_{\alpha_q} \in T_{\alpha_q} (T^\ast
Q)$, where $\pi_Q: T^\ast Q \rightarrow Q $ is the cotangent bundle
projection and $T_{ \alpha_q} \pi_Q: T_{ \alpha_q} (T^\ast Q) \rightarrow
T_q Q $ is its tangent map (or derivative) at $q$. In natural cotangent
bundle coordinates $(q^i, p_i)$, we have $\Theta_Q = p_i
\mathbf{d}q^i$ and $\Omega_Q = \mathbf{d}q^i\wedge \mathbf{d}p_i$.

Let $\Phi: G \times Q \rightarrow Q$ be a left smooth action of the Lie 
group $G$ on the manifold and $Q$. Denote by $g \cdot q = \Phi(g, q)$ the
action of $g \in G $ on the point $q \in Q $ and by $\Phi_g: Q
\rightarrow Q $ the diffeomorphism of $Q $ induced by $g $. The lifted
left action
$G
\times T^\ast Q \rightarrow T^\ast Q$, given by $g \cdot \alpha_q =
T^\ast_{g\cdot q}
\Phi_{g^{-1}} (\alpha_q)$  for $g \in G$ and $\alpha_q \in T^\ast _q Q$,
preserves $\Theta_Q$, and admits the equivariant momentum map
$\mathbf{J}: T^\ast  Q \rightarrow \mathfrak{g}^\ast$ whose expression is
$\langle \mathbf{J}(\alpha_q), \xi\rangle = \alpha_q((\xi_Q(q))$, where
$\xi \in \mathfrak{g}$, the Lie algebra of $G$, $\left\langle\,,
\right\rangle: \mathfrak{g}^\ast \times \mathfrak{g} \rightarrow
\mathbb{R}$ is the duality pairing between the dual $\mathfrak{g}^\ast$
and $\mathfrak{g}$, and $\xi_Q(q) = d\Phi(\exp t \xi, q)/dt|_{t=0}$ is
the value of the infinitesimal generator vector field $\xi_Q $ of the
$G$-action at $q \in Q$. Throughout this article it is assumed that  the
$G $-action on $Q$, and hence on $T^\ast Q$, is free and proper. Recall
also that $\left((T^\ast Q)_\mu, (\Omega_Q)_\mu\right)$ denotes the
reduced manifold at $\mu\in \mathfrak{g}^\ast$~\cite{or}, where $(T^\ast Q)_\mu : =
\mathbf{J}^{-1}( \mu) / G_ \mu$ is the orbit space of the $G_
\mu$-action on the momentum level manifold $\mathbf{J}^{-1}( \mu) $ and
$G_ \mu: = \{ g \in G \mid \operatorname{Ad}^\ast_g \mu = \mu \} $ is the
isotropy subgroup of the coadjoint representation of $G$ on
$\mathfrak{g}^\ast$. The left coadjoint representation of $g \in G $ on
$\mu\in \mathfrak{g}^\ast$ is denoted by $\operatorname{Ad}^\ast_{g^{-1}}
\mu$.

Cotangent bundle reduction at zero is already quite interesting
and has many applications. Let $\rho: Q \rightarrow Q/G $ be the
$G$-principal bundle projection defined by the proper free action of $G$
on $Q$, usually referred to as the {\bfi  shape space bundle}. Zero is a
regular value of
$\mathbf{J}$ and the map $\varphi_0:
\left((T^\ast Q)_0, (\Omega_Q)_0 \right) \rightarrow
\left(T^\ast (Q/G), \Omega_{Q/G} \right)$
given by  $\varphi_0([\alpha_q])(T_q \rho(v_q)): = \alpha_q(v_q)$, where
$\alpha_q \in \mathbf{J}^{-1}(0)$, $[\alpha_q] \in (T^\ast Q)_0$, and
$v_q \in T_q Q$, is a well-defined symplectic diffeomorphism.

This theorem generalizes in two non-trivial
ways when one reduces at a non-zero value of $\mathbf{J}$: an embedding
and a fibration  theorem. 
\medskip

\noindent \textbf{Embedding version of cotangent bundle reduction.}
Let $\mu \in \mathfrak{g}^\ast, Q_\mu := Q/G_\mu, \rho_
\mu: Q \rightarrow Q_\mu$ the projection onto the $G_
\mu$-orbit space, $\mathfrak{g}_\mu := \{\xi \in \mathfrak{g} \mid
\operatorname{ad}^\ast_\xi \mu = 0\}$ the Lie algebra of the coadjoint
isotropy subgroup $G_ \mu$, where $\operatorname{ad}_ \xi \eta: = [ \xi,
\eta] $ for any $\xi, \eta \in \mathfrak{g}$, $\operatorname{ad}^\ast_
\xi:\mathfrak{g}^\ast\rightarrow \mathfrak{g}^\ast$ the dual map, 
$\mu ^\prime := \mu |_{\mathfrak{g}_\mu}\in \mathfrak{g}_\mu
^\ast$  the restriction of $\mu$ to $\mathfrak{g}_\mu$, and $((T^\ast
Q)_\mu, (\Omega_Q) _\mu)$ the reduced space at $\mu$. The induced
$G_\mu$-action on $T^\ast Q$ admits the equivariant momentum map 
$\mathbf{J}^\mu : T^\ast Q \rightarrow \mathfrak{g}^\ast_\mu$ given by 
$\mathbf{J}^\mu(\alpha _q) = \mathbf{J} (\alpha _q) |_{\mathfrak{g}_\mu}$.
\textit{Assume there is a $G_\mu$-invariant one-form $\alpha _\mu$ on $Q$
with values in $(\mathbf{J}^\mu)^{-1}(\mu ^\prime)$}. Then there is a
unique closed two-form $\beta_\mu$ on $Q _\mu$ such that $\rho_\mu ^\ast 
\beta_\mu = \mathbf{d}\alpha_\mu$. Define the {\bfi magnetic term\/}  
$B_\mu := \pi^{\ast}_{Q_\mu} \beta _\mu$, where $\pi_{Q
_\mu} : T^\ast Q_\mu \rightarrow Q_\mu$ is the cotangent bundle
projection, which is a closed two-form on $T^\ast Q_ \mu$.
Then the map $\varphi_\mu : ((T^\ast Q)_\mu,(\Omega_Q) _\mu)
\rightarrow (T^\ast Q_\mu,\Omega_{Q_\mu}- B_\mu)$ given by
$\varphi_\mu([\alpha_q])(T_q\rho_\mu(v_q)) := (\alpha_q  -
\alpha_\mu(q))(v_q)$, for $\alpha_q \in \mathbf{J}^{-1}(\mu)$,
$[\alpha_q] \in (T ^\ast Q)_ \mu$, and
$v_q \in T_q Q$, is a symplectic embedding onto a submanifold of $T ^\ast
Q_ \mu$ covering the base $Q_ \mu$. The embedding $\varphi_\mu$ is a
diffeomorphism onto $T^\ast Q_\mu$ if and only if $\mathfrak{g} =
\mathfrak{g}_\mu$. If the one-form $\alpha_ \mu$ takes values in the
smaller set $\mathbf{J}^{-1}( \mu) $ then the image of $\varphi_ \mu$ is
the the vector subbundle $[T\rho_\mu(VQ)]^\circ$ of $T^\ast Q_{\mu}$,
where $VQ \subset TQ$ is the {\bfi vertical\/} vector subbundle consisting
of vectors tangent to the $G $-orbits, that is, its fiber at
$q \in Q $ equals $V_q Q = \{\xi_Q(q) \mid \xi \in \mathfrak{g}\}$, and
$^\circ $ denotes the annihilator relative to the natural duality pairing
between $T Q_{\mu} $ and $T ^{\ast} Q_{\mu}$. Note that if 
$\mathfrak{g}$ is Abelian or $\mu= 0$, the embedding $\varphi_\mu$ is
always onto and thus the reduced space is again, topologically, a
cotangent bundle.

It should be noted that there is a choice in this theorem, namely the one-form
$\alpha_\mu$. Whereas the reduced symplectic space $((T^\ast Q)_\mu,(\Omega_Q) _\mu)$
is intrinsic, the symplectic structure on the space $T^\ast Q_\mu$ depends on
$\alpha_\mu$. The theorem above states that no matter how $\alpha_\mu$ is
chosen, there is a symplectic diffeomorphism, which also depends on
$\alpha_\mu$, of the reduced space onto a submanifold of $T^\ast Q_\mu$. 
\medskip

\noindent \textbf{Connections.} The one-form $\alpha_ \mu$ is usually
obtained from a left connection on the principal bundle $\rho_\mu: Q
\rightarrow Q/G_\mu$ or $\rho: Q \rightarrow Q/G$. A left {\bfi
connection\/} one-form $A \in \Omega ^1(Q; \mathfrak{g})$ on the left
principal $G$-bundle $\rho: Q \rightarrow Q/G$ is a Lie algebra valued
one-form $A: TQ \rightarrow \mathfrak{g}$, where $\mathfrak{g}$ denotes
the Lie algebra of $G$, satisfying the conditions $A(\xi_Q) = \xi$ for
all $\xi\in \mathfrak{g}$ and $\mathcal{A}(T_q\Phi_g (v)) =
\operatorname{Ad}_g(\mathcal{A}(v))$ for all $g \in G $ and $v \in T_qQ$,
where $\operatorname{Ad}_g$ denotes the adjoint action of $G$ on
$\mathfrak{g}$. The {\bfi horizontal\/} vector subbundle $HQ$ of the
connection $A $ is defined as the kernel of $A $, that is, its fiber at
$q \in Q $ is the subspace $H_q:= \ker A(q)$. The map $v_q \mapsto
\operatorname{ver}_q (v_q):= [A(q) ( v _q ) ] _Q (q)$ 
is called the {\bfi vertical projection\/}, while the map $v_q \mapsto
\operatorname{hor}_q(v_q) := v_q - \operatorname{ver}_q(v_q)$ 
is called the {\bfi horizontal projection}. Since for any vector $v_q \in
T_q Q $ we have $v _q = \operatorname{ver}_q (v _q) +
\operatorname{hor}_q (v _q)$, it follows that $TQ= HQ \oplus VQ $
and the maps $ \operatorname{hor}_q : T_qQ \rightarrow H_qQ $ and
$\operatorname{ver}_q: T_qQ \rightarrow V_qQ$ are projections
onto the horizontal and vertical subspaces at every $q	 \in Q $.

Connections can be equivalently defined by the choice of a subbundle $HQ
\subset TQ$ complementary to the vertical subbundle $VQ $ satisfying
the following $G$-invariance property: $H_{g\cdot q}Q = T_q \Phi_g (H_qQ)$
for every $g \in G $ and $q \in Q $. The subbundle $H Q $ is called, as
before, the horizontal subbundle and a connection one-form $A $ is defined
by setting $A(q)(  \xi_Q(q) + u_q ) = \xi$, for any $\xi \in
\mathfrak{g}$ and $u_q \in H_qQ$. 

The {\bfi curvature\/} of the connection $A$ is the Lie algebra valued
two-form on $Q$ defined by $ 
B(u_q , v_q ) = \mathbf{d}A (\operatorname{hor}_q (u_q) ,
\operatorname{hor}_q (v _q ))$. When one replaces vectors in the exterior
derivative with their horizontal projections, then the result is called
the {\bfi exterior covariant derivative\/} and the preceding formula for
$B$ is often written as $B = DA$. Curvature {\it measures the lack of
integrability of the horizontal distribution\/}, namely, $B (u,v) = 
- A([\operatorname{hor} (u) , \operatorname{hor} (v)])$ for any two
vector fields $u$ and $v$ on $Q$. The {\bfi Cartan Structure Equations\/}
state that $B(u,v) = \mathbf{d}\mathcal{A} (u,v) - [\mathcal{A} (u) ,
\mathcal{A} (v) ]$, where the bracket on the right hand side is the Lie
bracket in $\mathfrak{g}$.

Since the connection $A$ is a Lie algebra valued one-form, for each
$\mu\in \mathfrak{g}^\ast$ the formula $\alpha_ \mu(q) := A(q)^\ast(\mu)$,
where $A(q)^\ast: \mathfrak{g} ^{\ast}\rightarrow T ^\ast _q Q$  is the
dual of the linear map $A(q) : T _q Q\rightarrow\mathfrak{g}$, defines a
usual one-form on $Q$. This one-form $\alpha_ \mu $ takes values in
$\mathbf{J} ^{-1}( \mu) $ and is equivariant in the following sense:
$\Phi _g ^\ast \alpha _\mu = \alpha _{ {\rm Ad } _g ^\ast \mu }$ for any
$g \in G$. 
\medskip

\noindent \textbf{Magnetic Terms and Curvature.} There are two methods to
construct the one-form $\alpha_ \mu$ from a connection. The first is to
start with a connection one-form $A^ \mu \in \Omega^1(Q; \mathfrak{g}_
\mu)$ on the principal $G_ \mu$-bundle $\rho_ \mu: Q \rightarrow Q/G_
\mu$. Then the one-form $\alpha_ \mu : = \left\langle \mu|_{
\mathfrak{g}_ \mu}, A ^\mu \right\rangle \in \Omega^1(Q) $ is $G_
\mu$-invariant and has values in $(\mathbf{J} ^\mu) ^{-1}( \mu|_{
\mathfrak{g}_ \mu})$. The magnetic term $B_ \mu$ is the pull back to
$T^\ast (Q/G_ \mu)$ of the $\mu|_{\mathfrak{g}_ \mu}$-component
$\mathbf{d}\alpha_\mu$ of the  curvature of $A^\mu$ thought of as a
two-form on the base $Q/G_ \mu$.

The second method is to start with a connection $A
\in \Omega^1(Q, \mathfrak{g})$ on the principal bundle $\rho: Q
\rightarrow Q/G $, to define $\alpha_ \mu : = \left\langle\mu, A
\right\rangle \in \Omega^1(Q) $, and to observe that this one-form is $G_
\mu$-invariant and has values  in $\mathbf{J}^{-1}( \mu) $. The magnetic 
term $B_ \mu$ is in this case the pull back to $T^\ast (Q/G_ \mu)$
of the $\mu$-component $\mathbf{d}\alpha_\mu$ of the  curvature of $A$
thought of as a two-form on the base $Q/G_ \mu$.

\medskip

\noindent \textbf{The Mechanical Connection.} If $(Q, \left\langle \!
\left\langle\,, \right\rangle \! \right\rangle)$ is a Riemannian manifold
and $G $ acts by isometries, there is a natural connection on the bundle
$\rho: Q \rightarrow  Q/G$, namely, define the horizontal  space at a
point to be the metric orthogonal to the vertical space. This
connection is called the {\bfi mechanical connection\/} and its
horizontal bundle consists of all vectors $v_q \in TQ$ such that
$\mathbf{J}( \left\langle \! \left\langle v_q, \cdot \right\rangle \!
\right\rangle) = 0 $. 

To determine the Lie algebra valued one-form $A $ of this connection, the
notion of {\bfi locked inertia tensor\/} needs to be introduced. This is
the  linear map $\mathbb{I}(q): \mathfrak{g} \rightarrow
\mathfrak{g}^\ast$ depending  smoothly on $q \in Q$ defined by the
identity $\left\langle
\mathbb{I}(q) \xi, \eta \right\rangle = \left\langle \! \left\langle
\xi_Q(q), \eta_Q(q) \right\rangle \! \right\rangle$ for any $\xi, \eta \in
\mathfrak{g}$. Since the $G $-action is  free, each $\mathbb{I}(q) $ is
invertible. The connection one-form whose horizontal space was
defined above is given by $A(q)(v_q) = \mathbb{I}(q)
^{-1}\left(\mathbf{J}( \left\langle \! \left\langle v_q, \cdot
\right\rangle \! \right\rangle ) \right)$.

Denote by $K: T ^\ast Q \rightarrow \mathbb{R}$ the kinetic  energy of
the metric $\left\langle \! \left\langle\,, \right\rangle \!
\right\rangle$ on the cotangent bundle, that is, $K( \left\langle \!
\left\langle v_q, \cdot \right\rangle \! \right\rangle) : =
\frac{1}{2}\|v_q \|^2 $. The one-form $\alpha_ \mu = A( \cdot )^\ast\mu$
is characterized for the mechanical connection $A$ by the condition 
${K}(\alpha _\mu (q)) = {\rm inf} \{ K(\beta_q ) \mid \beta_q \in
\mathbf{J}^{-1} (\mu ) \cap T_q ^\ast Q \}$.

\medskip

\noindent \textbf{The amended potential.} 
\label{the amended potiential}
A {\bfi simple mechanical system\/}
\index{mechanical!system!simple}%
\index{simple!mechanical system}%
is a Hamiltonian  system on a cotangent bundle $T ^\ast Q $ whose
Hamiltonian function is the the sum of the kinetic energy of a Riemannian
metric on $Q $ and a potential function  $V:Q \rightarrow \mathbb{R} $.
If there is a Lie group $G $ acting on $Q $ by isometries and leaving the
potential invariant, then we have a {\bfi simple mechanical system with
symmetry\/}.
\index{mechanical!system!simple with symmetry}%
\index{simple!mechanical system with symmetry}%
\index{symmetric!simple mechanical system}
The {\bfi amended\/} 
\index{amended potential}%
\index{potential!amended}%
or {\bfi effective potential\/}
\index{effective!potential}%
\index{potential!effective}%
$V_ \mu :Q \rightarrow \mathbb{R}$ at $\mu\in \mathfrak{g}^\ast$ is
defined by $V_\mu : = H \circ\alpha_\mu$, where
$\alpha_\mu$ is the one-form associated to the {\bfi mechanical
connection\/}.
\index{connection!mechanical}%
\index{mechanical!connection}%
Its expression in terms of the locked moment  of inertia tensor is given
by $V_{\mu }(q):=V(q)+\frac{1}{2}\langle
\mathbb{\mu },\mathbb{I}(q)^{-1}\mu \rangle$. The amended potential
naturally induces a smooth function $\widehat{V}_ \mu \in
C^\infty(Q/G_\mu)$.

The fundamental result about simple mechanical systems with
symmetry is the following. \textit{The push-forward by the embedding
$\varphi_\mu : ((T^\ast Q)_\mu,(\Omega_Q) _\mu)
\rightarrow (T^\ast Q_\mu,\Omega_{Q_\mu}- B_\mu)$ of the reduced 
Hamiltonian $H_\mu \in C^\infty((T^\ast Q)_\mu)$ of a simple
mechanical system $H =K + V \circ \pi_Q \in C ^\infty(T^\ast Q)$ is the
restriction to  the vector subbundle $\varphi_\mu((T^\ast Q)_\mu) \subset
T^\ast(Q/G_\mu)$, which is also a symplectic submanifold of
$(T^\ast(Q/G_\mu), \Omega_{Q/G_\mu} - B_\mu)$, of the simple mechanical
system on $T^\ast(Q/G_\mu)$ whose kinetic energy is given by the 
quotient Riemannian metric on $Q/G_\mu$ and whose potential is
$\widehat{V}_\mu$. However, Hamilton's equations on $T^\ast(Q/G_\mu)$ for
this simple mechanical system are computed relative to the magnetic
symplectic form $\Omega_{Q/G_\mu} - B_\mu$.}

There is a wealth of applications starting from this classical theorem to mechanical
systems, spanning such diverse areas as topological characterization of the level sets
of the energy-momentum map to methods of proving nonlinear stability of relative
equilibria (block diagonalization of the stability form in the application of
the energy-momentum method).
\medskip

\noindent \textbf{Fibration version of cotangent bundle reduction.} There
is a second theorem that realizes the reduced space of a cotangent bundle
as a locally trivial bundle over {\bfi shape space\/} $Q/G$. This version
is particularly well suited in the study of quantization problems and in
control theory. The result is the following. Assume that $G$ acts freely
and properly on $Q$. \textit{Then the reduced symplectic manifold $(T^\ast
Q)_\mu$ is a fiber bundle over $T^\ast (Q/G)$ with fiber the coadjoint
orbit $\mathcal{O}_\mu$}. How this is related to the Poisson structure of
the quotient $(T^\ast Q)/G $ will be discussed later.

\medskip

\noindent \textbf{The Kaluza-Klein construction.}
The extra term in the symplectic form of the reduced space is called a
magnetic term because it has this interpretation in electromagnetism. 
To understand why $B_\mu$ is called a
\textit{magnetic term\/}, consider the problem of a particle of mass $m$
and  charge $e$ moving in $\mathbb{R}^3$ under the influence of a given
magnetic field
$\mathbf{B} = B_x \mathbf{i} + B_y
\mathbf{j} + B_z
\mathbf{k}$,
\index{magnetic!field}
\index{field!magnetic}
$\operatorname{div} \mathbf{B} = 0$. The Lorentz Force Law (written in
the International System) gives the equations of motion
\begin{equation}
\label{Lorentz force law equations}
m \frac{d\mathbf{v}}{dt} = e\mathbf{v} \times \mathbf{B},
\end{equation}
where $e$ is the charge and $\mathbf{v} = (\dot{x}, \dot{y},
\dot{z}) =
\dot{\mathbf{q}}$ is the velocity of the particle. What is the Hamiltonian
description of these equations? 

There are two possible answers to this question.  To formulate them, 
associate to the divergence free vector field $\mathbf{B}$ the closed
two-form $B = B_x \mathbf{d}y\wedge \mathbf{d}z - B_y \mathbf{d}x \wedge
\mathbf{d}z + B_z \mathbf{d}x \wedge \mathbf{d}y$. Also, write
$\mathbf{B} = \operatorname{curl} \mathbf{A}$ for some other vector field
$\mathbf{A} = (A_x, A_y, A_z)$ on $ \mathbb{R}^3$, called the {\bfi
magnetic potential\/}.
\smallskip

\textit{Answer 1\/}. Take on $T^\ast \mathbb{R}^3$ the symplectic form 
$\Omega_B = \mathbf{d}x \wedge \mathbf{d}p_x + dy \wedge \mathbf{d}p_y +
\mathbf{d}z \wedge \mathbf{d}p_z - eB$, where
$(p_x, p_y, p_z ) = \mathbf{p} := m \mathbf{v}$ is the momentum of the
particle, and 
$h = m \| \mathbf{v} \|^2/2 = m(\dot{x}^2 + \dot{y}^2 + \dot{z}^2)/2$ is
the Hamiltonian, the kinetic  energy of the particle. A direct
verification shows that
$\mathbf{d}h = \Omega_B(X_h, \cdot)$, where
\begin{equation}
\label{Lorentz equations}
X_h = \dot{x} \frac{\partial}{\partial x} + 
\dot{y} \frac{\partial}{\partial y} +
\dot{z} \frac{\partial}{\partial z} 
+ e(B_z \dot{y} - B_y \dot{z})\frac{\partial}{\partial p_x}
+ e(B_x \dot{z} - B_z
\dot{x})\frac{\partial}{\partial p_y} + e(B_y \dot{x} - B_z
\dot{x})\frac{\partial}{\partial p_z},
\end{equation}
which gives the equations of motion \eqref{Lorentz force law equations}.
\smallskip

\textit{Answer 2\/}. Take on $T^\ast \mathbb{R}^3$ the canonical symplectic 
form $\Omega = \mathbf{d}x \wedge \mathbf{d}p_x + \mathbf{d}y \wedge
\mathbf{d}p_y + \mathbf{d}z \wedge \mathbf{d}p_z$ and the Hamiltonian
$h_\mathbf{A} = \|\mathbf{p} - e\mathbf{A}\|^2/2m$. A direct
verification shows that $\mathbf{d}h_\mathbf{A} = \Omega(X_{h_\mathbf{A}},
\cdot)$, where
$X_{h_\mathbf{A}}$ has the \textit{same expression\/} \eqref{Lorentz
equations}.

Next we show how the magnetic term in the symplectic form $\Omega_B$ is 
obtained by reduction from the {\bfi Kaluza-Klein\/} system. 
Let $Q = \mathbb{R}^3 \times S^1$ with the circle $G = S^1$ acting on $Q$,
only on the second factor. Identify the Lie algebra $\mathfrak{g}$ of
$S ^1$ with $\mathbb{R}$. Since the infinitesimal generator of this action
defined by $\xi \in \mathfrak{g} = \mathbb{R}$ has the
expression $\xi_Q(\mathbf{q}, \theta) = (\mathbf{q},\theta;
\mathbf{0},\xi)$, if $TS^1$ is trivialized as $S^1 \times \mathbb{R}$, a
momentum map
$\mathbf{J}: T^\ast Q = \mathbb{R}^3 \times S^1 \times \mathbb{R}^3
\times\mathbb{R} \rightarrow \mathfrak{g}^\ast = \mathbb{R}$ is given by
$\mathbf{J}(\mathbf{q},\theta; \mathbf{p}, p) \xi = 
(\mathbf{p}, p) \cdot (\mathbf{0}, \xi) = p \xi$, that is,
$\mathbf{J}(\mathbf{q}, \theta; \mathbf{p}, p) = p$. In this case the
coadjoint action is trivial, so for any $\mu \in \mathfrak{g}^\ast =
\mathbb{R}$, we have $G_\mu = S^1$, $\mathfrak{g}_\mu = \mathbb{R}$, 
and $\mu' = \mu$. The one-form $\alpha_\mu = \mu (A_x \mathbf{d}x + A_y
\mathbf{d}y + A_z \mathbf{d}z + \mathbf{d}\theta) \in \Omega^1(Q)$, where
$d\theta$ denotes the length one-form on $S^1$, is clearly $G_\mu =
S^1$-invariant, has values in $\mathbf{J}^{-1}(\mu) = \{(\mathbf{q},
\theta; \mathbf{p}, \mu) \mid \mathbf{q}, \mathbf{p} \in \mathbb{R}^3,
\theta \in S^1 \}$, and its exterior differential equals
$\mathbf{d}\alpha_\mu = \mu B$. Thus the closed two-form
$\beta_\mu$ on the base $Q_\mu = Q/G_\mu = Q/S^1= \mathbb{R}^3$ equals 
$\mu B$ and hence the magnetic term, that is, the closed two-form 
$B_\mu = \pi_{Q_\mu}^\ast \beta_\mu$ on $T^\ast Q_ \mu =
T^\ast\mathbb{R}^3 $, is also $\mu B$ since
$\pi_{Q_\mu}: Q = \mathbb{R}^3 \times S^1 \rightarrow Q/G_\mu =
\mathbb{R}^3$ is the projection. Therefore, the reduced space
$(T^\ast Q)_\mu$ is symplectically diffeomorphic to 
$(T^\ast \mathbb{R}^3, \mathbf{d}x\wedge \mathbf{d}p_x+ \mathbf{d}y
\wedge \mathbf{d}p_y + \mathbf{d}z \wedge \mathbf{d}p_z -
\mu B)$, which coincides with the phase space in
\textit{Answer 1\/} if we put $\mu = e$. This also gives the
physical interpretation of the momentum map $\mathbf{J}: T^\ast Q = 
\mathbb{R}^3 \times S^1 \times \mathbb{R}^3
\times\mathbb{R} \rightarrow \mathfrak{g}^\ast = \mathbb{R}$,
$\mathbf{J}(\mathbf{q}, \theta; \mathbf{p}, p) = p$ and hence of the
variable conjugate to the circle variable $\theta$: $p$ \textit{represents
the charge}. Moreover, \textit{the magnetic term in
the symplectic form is, up to a charge factor, the magnetic field.\/} 

The kinetic energy Hamiltonian $h(\mathbf{q}, \theta; \mathbf{p}, p) :=
\frac{1}{2m}\|\mathbf{p}\|^2 + \frac{1}{2} p^2$ of the Kaluza-Klein
metric, that is, the Riemannian metric obtained by keeping the standard
metrics on each factor and declaring  $\mathbb{R}^3$ and $S^1$
orthogonal, induces the reduced Hamiltonian $h_\mu(\mathbf{q}) =
\frac{1}{2m}\|\mathbf{p}\|^2 + \frac{1}{2}\mu^2$ which, up to the
constant $\mu^2/2$, equals the kinetic energy Hamiltonian in
\textit{Answer 1\/}. Note that this reduced system is \textit{not\/} the
geodesic flow of the Euclidean metric because of the presence of the
magnetic term in the symplectic form.  However, \textit{the equations of
motion of a charged particle in a magnetic field are obtained by reducing
the geodesic flow of the Kaluza--Klein metric.}

A similar construction is carried out in Yang--Mills theory where $A$ is a connection on
a principal bundle and $B$ is its curvature. Magnetic terms appear also in classical
mechanics. For example, in rotating systems the Coriolis force (up to a dimensional factor)
plays the role of the magnetic term.

\medskip

\noindent \textbf{Reconstruction of dynamics for cotangent bundles.} 
A general reconstruction method  of the dynamics from the reduced
dynamics was given in \cite{or}. For cotangent bundles, using the
mechanical connection, this method simplifies considerably. 

Start with the following general situation. Let $G $ act freely on the
configuration manifold $Q$,   $h: T^\ast Q
\rightarrow \mathbb{R}$ be  a $G$-invariant Hamiltonian, $\mu \in
\mathfrak{g}^\ast$, $\alpha_q \in \mathbf{J}^{-1}(\mu)$, and $c_\mu(t)$  the
integral curve of the reduced system  with initial condition $[\alpha_q]
\in (T^\ast Q)_\mu$ given by the reduced Hamiltonian function $h_\mu:
(T^\ast Q)_\mu
\rightarrow \mathbb{R}$. In terms of a connection $A \in
\Omega^1(\mathbf{J}^{-1}(\mu); \mathfrak{g}_\mu)$ on the left
$G_\mu $-principal bundle $\mathbf{J}^{-1}(\mu) \rightarrow (T^\ast Q)_\mu$
the reconstruction procedure proceeds in four steps:
\begin{itemize}
\item \textbf{Step 1:} \textit{Horizontally lift the curve $c_\mu(t) \in (T^\ast
Q)_\mu$ to a curve $d(t) \in \mathbf{J}^{-1}(\mu)$ with $d(0)= \alpha_q $;}
\item \textbf{Step 2:} \textit{Set $\xi(t) = A(d(t))\left(X_h(d(t))\right)
\in
\mathfrak{g}_\mu$;}
\item \textbf{Step 3:} \textit{With $\xi(t) \in \mathfrak{g}_\mu $
determined in  Step 2, solve the nonautonomous differential equation
$\dot{g}(t) = T_e L_{g(t)}
\xi(t)$ with initial condition $g(0) = e$, where $L_ g $ denotes left
translation on $G $; this is the step that  involves ``quadratures" and
is the main obstacle to finding  explicit formulas;}
\item \textbf{Step 4:} \textit{The curve $c(t) = g(t)\cdot d(t)$, with
$d(t)$ found in Step 1 and $g(t)$ found in Step 3 is the integral curve
of $X_h$ with initial condition $c(0) = \alpha_q $.}
\end{itemize}
This method depends on the choice of the connection $A \in
\Omega^1(\mathbf{J}^{-1}(\mu); \mathfrak{g}_\mu)$. Here are several
particular cases when this procedure simplifies.
\smallskip

\noindent \textbf{(a)} \textit{One-dimensional coadjoint isotropy group.}
\label{one dimensional coadjoint isotropy group}
If $G_\mu = S^1$ or $G_\mu = \mathbb{R}$, identify $\mathfrak{g}_ \mu $
with $\mathbb{R}$ via the map $a \in \mathbb{R}\leftrightarrow a \zeta \in
\mathfrak{g}_ \mu$, where $\zeta \in \mathfrak{g}_\mu$, $\zeta\neq 0$, is
a generator of $\mathfrak{g}_ \mu$. Then a connection one-form on the
$S^1$ (or $\mathbb{R}$) principal bundle $\mathbf{J}^{-1}(\mu)
\rightarrow (T^\ast Q)_\mu$ is the  one-form $A$ on
$\mathbf{J}^{-1}(\mu)$  given by $A = \frac{1}{\langle \mu, \zeta\rangle}
\theta_\mu$, where $ \theta_ \mu$ is the pull-back of the canonical
one-form $\theta \in\Omega^1(T^\ast Q)$ to the submanifold
$\mathbf{J}^{-1}(\mu)$. The curvature of this connection is the two-form
on $( T^\ast Q) _\mu$ given by $ \operatorname{curv}(A) = -
\frac{1}{\langle \mu, \zeta\rangle} \omega_\mu$,
where $\omega_ \mu$ is the reduced symplectic form on $(T^\ast Q)_\mu $.
In this case, the curve $\xi(t) \in \mathfrak{g}_\mu$ in Step 2 is given
by $ \xi(t) = \Lambda[h](d(t))$,
where $\Lambda \in \mathfrak{X}(T^\ast Q)$ is the {\bfi Liouville vector
field\/}
\index{vector!field!Liouville}%
\index{Liouville vector field}%
characterized by the property of being the unique vector field on $T^\ast Q$
that satisfies  the relation $\mathbf{d}\theta( \Lambda, \cdot ) = \theta $.
In canonical coordinates $(q^i, p_i)$ on $T^\ast Q$, $\Lambda = p_i
\frac{\partial}{\partial p_i}$.

\noindent \textbf{(b)} \textit{Induced connection.}
\label{induced connection}
Any connection $\mathcal{A} \in \Omega^1(Q; \mathfrak{g}_\mu)$ on the left principal
bundle $Q \rightarrow Q/G_\mu$ induces a connection $A \in
\Omega^1(\mathbf{J}^{-1}(\mu); \mathfrak{g}_\mu)$ by
$A(\alpha_q)\left(V_{\alpha_q}\right) : = \mathcal{A}(q) \left(
T_{\alpha_q} \pi_ Q\left(V_{\alpha_q}\right)\right)$,
where $q \in Q$, $\alpha_ q \in T_q^\ast Q$, $V_{\alpha_q} \in T_{\alpha_q}(T^\ast
Q)$, and $\pi_Q: T^\ast Q \rightarrow  Q $ is the cotangent bundle projection. In
this case,  the curve $\xi(t) \in \mathfrak{g}_\mu$ in Step 2 is given by 
$\xi(t)= \mathcal{A}(q(t))\big( \mathbb{F}h(d(t)\big)$,
where $q(t): = \pi_Q(d(t))$ is the base integral curve and the vector
bundle morphism $\mathbb{F}h: T^\ast Q \rightarrow TQ$ is the fiber derivative of $h$
given by $
\mathbb{F}h(\alpha_q)(\beta_q): = \left.\frac{d}{dt}\right|_{t=0}h(\alpha_q +
t\beta_q)$
for any $\alpha_q, \beta_ q \in T_a^\ast Q $.
\smallskip
Two particular instances of this situation are noteworthy.
\smallbreak

\textbf{(b1)} Assume that the Hamiltonian $h$ is that of a
simple mechanical system with symmetry. Choosing 
$A$ to be the mechanical connection $A_{\operatorname{mech}}$,  the curve
$\xi(t) \in \mathfrak{g}_\mu$ in Step 2 is given by $\xi(t)=
A_{\operatorname{mech}}(q(t))\left(\langle\!\langle d(t),
\cdot \rangle \!\rangle\right)$. 
\smallbreak

\textbf{(b2)} If $Q = G$ is a Lie group,
$\operatorname{dim}G_\mu = 1$, and
$\zeta$ is a generator of $\mathfrak{g}_\mu$, then the connection $A
\in \Omega ^1(G)$ can be chosen to equal $\mathcal{A}(g): =
\frac{1}{\langle \mu, \zeta\rangle} T^\ast _g R_{g^{-1}}(\mu)$, where
$\zeta$ is a generator of $\mathfrak{g}_ \mu$ and $R_g $is right
translation on  $G $.

\noindent \textbf{(c)} \textit{Reconstruction of dynamics for simple
mechanical systems with symmetry.}
The case of simple mechanical systems with symmetry deserves special 
attention since several steps in the reconstruction method can be
simplified. For simple mechanical systems the knowledge of the base
integral curve $q(t) $ suffices to determine the entire integral curve on
$T^\ast Q$. Indeed, if $h = K + V\circ \pi_q$ is the Hamiltonian, the
Legendre transformation $\mathbb{F}h: T^\ast Q \rightarrow TQ $ determines
the Lagrangian system on $TQ$ given by $\ell(u_q) = \frac{1}{2}\|u_q\|^2
- V(u_q)$, for $u_q
\in T_q Q$. Lagrange's equations are second order and thus the evolution
of the velocities  is given by the time derivative $\dot{q}(t)$ of the
base  integral curve. Since $\mathbb{F}h =
(\mathbb{F}\ell)^{-1}$, the solution of the Hamiltonian system
is given by $\mathbb{F}\ell (\dot{q}(t))$. Using the explicit expression
of the mechanical connection and the notation given in the general
procedure, the method of reconstruction simplifies to the following steps.
To find the integral curve $c(t)$ of the simple mechanical system with
$G$-symmetry $h = K + V \circ \pi_Q$ on $T^\ast Q$ with initial condition
$c(0) = \alpha_q \in T^\ast _q Q$, knowing the integral curve $c_\mu(t)$
of the reduced Hamiltonian system on $(T^\ast Q)_\mu$ given by the
reduced Hamiltonian function $h_\mu : (T^\ast Q)_\mu 
\rightarrow \mathbb{R}$ with initial condition $c_\mu(0) = [\alpha_q]$ one proceeds
in the following manner. Recall the symplectic embedding 
$\varphi_\mu : \bigl((T^\ast Q)_\mu,\linebreak (\Omega_Q) _\mu\bigr) 
\rightarrow (T^\ast (Q/G_\mu), \Omega_{Q/G_\mu}- B_\mu)$. The curve $\varphi_\mu (c_\mu(t))
\in T^\ast(Q/G_\mu)$ is an integral curve of the Hamiltonian system on
$\left(T^\ast (Q/G_\mu), \Omega_{Q/G_\mu} - B_\mu\right)$ given by the function that
is the sum of the kinetic energy of the quotient Riemannian metric and the quotient
amended potential $\widehat{V}_\mu$. Let $q_\mu(t) : =
\pi_{Q/G_\mu}(c_\mu(t))$ be the base integral curve of this system, where
$\pi_{Q/G_\mu}: T^\ast(Q/G_\mu) \rightarrow Q/G_\mu$ is the cotangent
bundle projection.

\begin{itemize}

\item \textbf{Step 1.} \textit{Relative to the mechanical connection
$A_{\operatorname{mech}}
\in \Omega^1 (Q; \mathfrak{g}_\mu)$, horizontally lift  $q_\mu(t) \in Q/G_\mu$ to a
curve $q_h(t) \in Q$ passing through $q_h(0) = q$.} 

\item \textbf{Step 2.} \textit{Determine $\xi(t)\! \in\! \mathfrak{g}_\mu
$ from the algebraic system $\!\langle\!\langle \xi(t)_Q(q_h(t)),\!
\eta_Q(q_h(t))
\!\rangle\!\rangle\linebreak
 = \langle \mu, \eta \rangle$ for all $\eta \in
\mathfrak{g}_\mu
$, where $\langle\!\langle \cdot, \cdot \rangle\!\rangle$ is the $G$-invariant
kinetic energy Riemannian metric on $Q$. This implies that $\dot{q}_h (0)$
and
$\xi(0)_Q(q)$ are the horizontal and vertical components of the vector $\alpha_q
^\sharp \in T_q Q$ which is associated by the metric $\langle\!\langle \cdot, \cdot
\rangle\!\rangle$ to the initial condition $\alpha_q $.}

\item \textbf{Step 3.} \textit{Solve 
$\dot{g}(t) = T_e L_{g(t)} \xi(t)$ in $G_\mu $ with initial condition
$g(0) = e $.}

\item \textbf{Step 4.} \textit{The curve $q(t) := g(t) \cdot q_h(t)$, with
$q_h(t) $ and $g(t)$ determined in Steps 2 and 4 respectively, is the
base integral curve of the simple mechanical system with symmetry defined
by the function $h $ satisfying
$q(0) = 0 $. The curve $(\mathbb{F}h)^{-1} (\dot{q}(t)) \in T^\ast Q$ is
the integral curve of this system with initial condition
$c(0) = \alpha_q$. In addition, $q'(t) = g(t) \cdot \big(\dot{q}_h(t) +
\xi(t)_Q(q_h(t))\big)$ is the horizontal plus vertical decomposition relative to
the connection induced on $\mathbf{J}^{-1}(\mu) \rightarrow (T^\ast Q)_\mu$ by the
mechanical connection $A_{\operatorname{mech}} \in \Omega^1(Q;
\mathfrak{g}_\mu)$. } 
\end{itemize} 

There are several important situations when Step 3, the main obstruction
to an explicit solution of the reconstruction problem, can be carried
out. We shall review some of them below.  
 
\textbf{(c1)} \textit{The case $G_\mu = S^1$.}
If $G_\mu$ is Abelian, the equation in Step 3 has the solution $g(t) =
\exp \int_0^t \xi(s) ds$. If, in addition, $G_\mu = S^1$, then $\xi(s)$
can be explicitly determined by Step 2. Indeed, if $\zeta \in
\mathfrak{g}_\mu $ is a generator of
$\mathfrak{g}_\mu $, writing $\xi(s) = a(s) \zeta$ for some smooth real valued
function $a$ defined on some open interval around  the origin, the algebraic equation
in Step 2 implies that $\langle\!\langle a(s)\xi(t)_Q(q_h(t)),
\zeta_Q(q_h(t)) \rangle\!\rangle = \langle \mu, \zeta \rangle$ which gives
$ a(s) = \langle \mu, \zeta\rangle/\| \zeta_Q(q_h(s))\|^2$.
Therefore, the base integral curve of the solution of the simple mechanical system
with symmetry on $T^\ast Q$ passing through $q $ is
\[
q(t) = \exp \left( \langle \mu, \zeta\rangle \int_0^t \frac{ ds}{\|
\zeta_Q(q_h(s))\|^2} \zeta \right) \cdot q_h(t)
\]
and 
\[
\dot{q}(t) = \exp \left( \langle \mu, \zeta\rangle \int_0^t \frac{ ds}{\|
\zeta_Q(q_h(s))\|^2} \zeta \right) \cdot \left(\dot{q}_h(t)  + 
\frac{\langle \mu, \zeta\rangle}{\| \zeta_Q(q_h(s))\|^2} \zeta_Q(q_h(t)) \right).
\]

\textbf{(c2)} \textit{The case of compact Lie groups.}
An obvious situation when the differential equation in Step 3 can be
solved is if $\xi(t) = \xi$ for all $t$, where $\xi$ is a given element
of $\mathfrak{g}_\mu$. Then the solution is $g(t) = \exp (t \xi)$.
However, Step 2 puts certain restrictions under this hypothesis, because
it requires that $\langle\!\langle \xi(t)_Q(q_h(t)),
\eta_Q(q_h(t)) \rangle\!\rangle = \langle \mu, \eta \rangle$ for any $\eta \in
\mathfrak{g}_\mu $. This is satisfied if there is a bilinear nondegenerate form
$(\cdot , \cdot )$ on $\mathfrak{g}$ satisfying $(\zeta, \eta ) =
\langle\!\langle \zeta_Q(q),\eta_Q(q) \rangle\!\rangle$ for all $q \in Q$
and $\zeta, \eta \in \mathfrak{g}$.
This implies that $(\cdot , \cdot )$ is positive
definite and invariant under the adjoint action of $G$ on $\mathfrak{g}$
so semisimple Lie algebras of noncompact type are excluded. If
$G$ is compact, which ensures the existence of a positive adjoint
invariant inner product on $\mathfrak{g}$, and $Q = G$, this condition
implies that the kinetic energy metric is invariant under the adjoint
action. There are examples in which such conditions are
natural, such as in Kaluza--Klein theories.
Concluding, if $G$ is a compact Lie group and $(\cdot, \cdot)$ is a
positive definite metric invariant under the adjoint action of $G$ on
$\mathfrak{g}$ satisfying $(\zeta, \eta ) =
\langle\!\langle \zeta_Q(q),\eta_Q(q) \rangle\!\rangle$ for all $q \in Q$
and $\zeta, \eta \in \mathfrak{g}$, then the element $\xi(t)$
in Step 2 can be chosen to be constant and is determined by the identity
$(\xi, \cdot) = \mu|_{\mathfrak{g}_\mu}$ on $\mathfrak{g}_\mu$. The
solution of the equation on Step 3 is then $g(t) = \exp (t\xi)$. 

\textbf{(c3)} \textit{The case when $\dot{\xi}(t)$ is proportional to
$\xi(t)$.}
Try to find a real valued function $f(t) $ such that $g(t) = \exp(f(t) \xi(t)) $ is a
solution of the equation $\dot{g}(t) = T_e L_{g(t)} \xi(t) $ with $f(0) =
0 $. This gives, for small $t $, the equation $\dot{f}(t)
\xi(t) + f(t) \dot{\xi}(t) = \xi(t) $, that is, it is necessary that  
$\xi(t)$ and $\dot{\xi}(t)$ be proportional. So if $\dot{\xi}(t) =
\alpha(t) \xi(t)$ for some known smooth function $\alpha(t) $, then this
gives $f(t) = \int_0^t \exp \left( \int_t^s \alpha(r) dr \right) ds$.

\textbf{(c4)} \textit{The case of $G_\mu $ solvable.}
Write $g(t) = \exp(f_1(t) \xi_1) \exp(f_2(t) \xi_2) \dots \exp(f_n(t)
\xi_n)$, for some basis $\{ \xi_1, \xi_2, \dots, \xi_n \} $ of
$\mathfrak{g}_\mu $ and some smooth real valued functions $f_i $, $i =
1,2,\dots , n $, defined around zero. It is known that if $G_\mu $ is
solvable, the equation in Step 3 can be solved by quadratures for the
$f_i$. 

\medskip

\noindent \textbf{Reconstruction phases for simple mechanical systems with
$S^1 $ symmetry.} Consider a simple mechanical system with symmetry $G $
on the Riemannian manifold $(Q, \langle\!\langle \cdot, \cdot
\rangle\!\rangle)$ with $G $-invariant potential $V \in C^\infty(Q)$. If
$\mu\in \mathfrak{g}^\ast$ let $V_ \mu$ be the amended
potential and $\widehat{V}_ \mu\in C^\infty(Q/G_\mu)$ the induced function
on the base. Let $c: [0, T] \rightarrow T ^\ast Q $ be an integral
curve of the system with Hamiltonian $h = K + V \circ \pi_Q$ and suppose
that its projection $c_\mu: [0, T] \rightarrow (T ^\ast Q)_\mu$ to the
reduced space is a closed integral curve of the reduced system with
Hamiltonian $h_\mu$. The {\bfi reconstruction phase\/}
\index{reconstruction!phase}%
\index{phase!reconstruction}%
associated to the loop $c_\mu(t)$ is the group element $g \in G_\mu$, satisfying the
identity $c(T)= g \cdot c(0)$. We shall present two explicit formulas of
the reconstruction phase for the case when $G_\mu = S^1$. Let
$\zeta \in \mathfrak{g}_\mu = \mathbb{R}$ be a generator of the coadjoint isotropy
algebra and write $c(T) = \exp(\varphi \zeta) \cdot c(0)$; in this case $\varphi$ is
identified with the reconstruction phase and, as we shall see in concrete
mechanical examples, it truly represents an angle.

If $G_\mu = S^1 $, the $G_\mu$-principal bundle $\pi_\mu: \mathbf{J}^{-1}(\mu)
\rightarrow (T^\ast Q)_\mu: =\mathbf{J}^{-1}(\mu)/G_\mu$ admits two natural
connections: $A =\frac{1}{\mu \zeta}\theta_\mu \in
\Omega^1(\mathbf{J}^{-1}(\mu))$,  where $\theta_\mu $ is the pull-back of
the canonical one-form on the cotangent bundle to the momentum level
submanifold $\mathbf{J}^{-1}(\mu)$,  and $\pi_Q^\ast
A_{\operatorname{mech}} \in \Omega^1(\mathbf{J}^{-1}(\mu))$.
There is no reason to choose one connection over the other and thus there
are two natural formulas for the reconstruction phase in this case. Let 
$c_\mu(t) $ be a periodic orbit of period $T $ of the reduced
system and denote also by $h_\mu$ the value of the Hamiltonian function on it. Assume 
that $D$ is a two-dimensional surface in $(T^\ast Q)_\mu$ whose boundary is the loop
$c_\mu(t)$. Since the manifolds $(T ^\ast Q)_\mu$ and $T ^\ast (Q/S^1)$
are diffeomorphic (but not symplectomorphic), it makes sense to consider
the base integral curve $q_\mu(t)$ obtained by projecting $c_\mu(t) $ to
the base $Q/S^1 $, which is a closed curve of period $T $. Denote by 
$ \langle \widehat{V}_\mu\rangle := \frac{1}{T}\int_0^T
\widehat{V}_\mu(q_\mu(t)) dt$  the average of $\widehat{V}_\mu$ over the
loop $q _\mu(t)$. Let $q_h(t) \in Q$ be the
$A_{\operatorname{mech}}$-horizontal lift of $q_\mu(t) $ to $Q$
and let $\chi $ be the $A_{\operatorname{mech}}$-holonomy of the
loop $q_\mu(t)$ measured from $q(0)$, the base point of $c(0)$; its
expression is given by $\exp \chi = \exp\left(- \int\!\!\!\int_D
B\right)$, where $B $ is the curvature of the mechanical connection.
Denote by $\omega_ \mu$ the reduced symplectic form on $(T^\ast Q)_\mu$.
With these notations the phase $\varphi$ is given by:
\begin{equation}
\label{general phase formula}
\varphi = \frac{1}{\mu\zeta}\int \!\!\!\int_D \omega_\mu + \frac{2(h_\mu
- \langle \widehat{V}_\mu \rangle)T}{\mu\zeta} 
= \chi + \mu\zeta\int_0^T \frac{ds}{\|\zeta_Q(q_h(s))\|^2}\,.
\end{equation}
The first terms in both formulas are the so-called {\bfi geometric phases\/}
\index{geometric phase}%
\index{phase!geometric}%
because they carry only geometric information given by the connection, whereas the
second terms are called the {\bfi dynamic phases\/}
\index{dynamic phase}%
\index{phase!dynamic}%
since they encapsulate information directly linked to the Hamiltonian. The expression of
the total phase as a sum of a geometric 
 and a dynamic phase is not intrinsic and is
connection dependent. It can even happen that one of these summands vanishes.
We shall consider now two concrete examples: the free rigid body and the
heavy top.

\medskip
\noindent \textbf{Reconstruction phases for the free rigid body.}
\label{phases for the free rigid body}
The motion of the free rigid body is a geodesic with respect to a left invariant
Riemannian metric on $\operatorname{SO}(3)$ given by the moment of inertia of the
body. The
phase space of the free rigid body motion is $T ^\ast\operatorname{SO}(3)$ and a
momentum map $\mathbf{J}:T ^\ast\operatorname{SO}(3) \rightarrow \mathbb{R}^3$ of the
lift of left translation to the cotangent bundle is given by right translation to the
identity element. We have identified here $\mathfrak{so}(3)$ with
$\mathbb{R}^3$ by the Lie algebra isomorphism $\mathbf{x} \in (\mathbb{R}^3, \times)
\mapsto \hat{\mathbf{x}} \in (\mathfrak{so}(3), [\cdot, \cdot ])$, where
$\hat{\mathbf{x}}(\mathbf{y}) = \mathbf{x} \times \mathbf{y}$, and
$\mathfrak{so}(3)^\ast$ with $\mathbb{R}^3$ by the inner product on $\mathbb{R}^3$.
The reduced manifold $\mathbf{J}^{-1}(\mu)/G_\mu$ is identified with the sphere
$S^2_{\|\mu \|}$ in $\mathbb{R}^3$ of radius $\| \mu\| $ with the symplectic form
$\omega_\mu = -dS/\|\mu\|$, where $dS$ is the standard area form on $S^2_{\|\mu \|}$
and where $G_\mu \cong S^1 $ is the group of rotations around the axis $\mu$. These
concentric  spheres are the coadjoint orbits of the Lie--Poisson space
$\mathfrak{so}(3)^\ast$ and represent the level sets of the Casimir functions that
are all smooth functions of $\| \Pi\|^2 $, where $\Pi\in
\mathbb{R}^3 $ denotes the body angular momentum.

The Hamiltonian of the rigid body on the Lie--Poisson space $T ^\ast
\operatorname{SO}(3)/\operatorname{SO}(3) \cong \mathbb{R}^3 $ is given by 
$ h(\Pi):= \frac{1}{2}\left(\frac{\Pi_1 ^2}{I_1}  + \frac{\Pi_2 ^2}{I_2} +
\frac{\Pi_3 ^2}{I_3} \right) $ where $I_1, I_2, I_3 >0 $ are the principal
moments of inertia of the body. Let
$\mathbb {I} : = \operatorname{diag}(I_1, I_2, I_3)$ denote the
moment of inertia  tensor diagonalized
in a principal axis body frame. The Lie--Poisson bracket on
$\mathbb{R}^3$ is given by $\{f, g\}(\Pi) = - \Pi\cdot \left( \nabla
f(\Pi) \times \nabla g(\Pi) \right)$ and the equation of motions are
$\dot{\Pi} = \Pi\times \Omega$,  where $\Omega \in \mathbb{R}^3$ is the
body angular velocity given in terms of $\Pi$ by $\Omega_i := \Pi/I_i $,
for $i = 1,2,3 $, that is, $\Omega= \mathbb {I}^{-1} \Pi$. The
trajectories of the these equations are found by intersecting a family of
homothetic energy ellipsoids with the angular momentum concentric
spheres. If $I_1> I_2> I_3 $, one immediately sees that all orbits are
periodic with the exception of four centers (the two possible rotations
about the long and the  short moment of inertia axis of the body), two
saddles (the two rotations about the middle moment of inertia axis of the
body), and four heteroclinic orbits connecting the two saddles.

Suppose that $\Pi(t)$ is a periodic orbit on the sphere $S^2_{\|\mu\|}$ with period
$T $. After time $T$ by how much has the rigid body rotated in space? The answer to
this question follows directly from \eqref{general phase formula}. Taking $\zeta =
\mu/\| \mu\| $ and the potential $v \equiv 0 $ we get 
\begin{align*}
\varphi&= - \Lambda + \frac{2h_\mu T}{\| \mu\|}
=\int\!\!\!\int_D
\frac{2\| \mathbb {I} \Pi(s)\|^2 - (\Pi(s) \cdot \mathbb
{I}\Pi(s))(\operatorname{tr}
\mathbb {I}) }{ (\Pi(s) \cdot \mathbb {I}\Pi(s))^2}ds + \| \mu\|^3 \int_0^T
\frac{ds}{(\Pi(s) \cdot \mathbb {I}\Pi(s))},
\end{align*} 
where $D$ is one of the two spherical caps on $S^2_{\| \mu\|} $ whose boundary 
is the periodic orbit $\Pi(t)$, $h_\mu$ is the value of the total energy on the
solution $\Pi(t)$, and $\Lambda$ is the oriented solid angle, that is, 
\[
\Lambda: = - \frac{1}{\| \mu\|} \int \!\!\!\int_D \omega_\mu, \qquad | \Lambda | =
\frac{\operatorname{area} D}{\| \mu\|^2}.
\]

\noindent \textbf{Reconstruction phases for the heavy top.}
The heavy top is a simple mechanical systems with symmetry $S^1$ on
$T ^\ast\operatorname{SO}(3) $ whose Hamiltonian function is given by
$h(\alpha_h) : = \frac{1}{2}\| \alpha_h  ^\sharp \|^2 + Mg \ell
\mathbf{k} \cdot h\chi$, 
where $h \in\operatorname{SO}(3)$, $\alpha_h \in T ^\ast_h \operatorname{SO}(3)$, 
$\mathbf{k}$ is the unit vector of the spatial $Oz$ axis (pointing in the opposite
direction of the gravity force), $M \in \mathbb{R}$ is the total mass of the body, $g
\in \mathbb{R}$ is the value of the gravitational acceleration, the fixed point about
which the body moves is the origin, and $\chi $ is the unit vector of the
straight line segment of length $\ell $ connecting the origin to the center of mass
of the body. This Hamiltonian is left invariant under rotations about the spatial
$Oz$ axis. A momentum map induced by this $S^1$-action is given by $\mathbf{J}: T
^\ast\operatorname{SO}(3) \rightarrow \mathbb{R}$, $\mathbf{J}(\alpha_h) = 
T ^\ast_e L_h(\alpha_h) \cdot \mathbf{k}$; recall that $T^\ast_e
L_h(\alpha_h) = : \Pi
\in \mathbb{R}^3$ is the body angular momentum. The reduced space
$\mathbf{J}^{-1}(\mu)/S^1$ is generically the cotangent bundle of the unit
sphere endowed with the symplectic structure given by the
sum of the canonical form plus a magnetic term; equivalently, this is the coadjoint orbit in the dual of the Euclidean Lie
algebra
$\mathfrak{se}(3)^\ast = \mathbb{R}^3 \times \mathbb{R}^3$ given by $\mathcal{O}_\mu=
\{(\Pi, \Gamma) \mid \Pi\cdot \Gamma = \mu, \| \Gamma \|^2 = 1 \} $. The projection
map $\mathbf{J}^{-1}(\mu) \rightarrow \mathcal{O}_\mu$ implementing the symplectic
diffeomorphism between the reduced space and the coadjoint orbit in
$\mathfrak{se}(3)^\ast $ is given by $\alpha_h \mapsto (\Pi, \Gamma) :=(T_e ^\ast L_h
(\alpha_h), h ^{-1}\mathbf{k})$. The orbit symplectic form $\omega_\mu $ on
$\mathcal{O}_\mu$ has the expression $\omega_\mu(\Pi, \Gamma)((\Pi\times
\mathbf{x} + \Gamma \times \mathbf{y},
\Gamma\times \mathbf{x}), (\Pi\times \mathbf{x}' + \Gamma \times \mathbf{y}',
\Gamma\times \mathbf{x}')) = -\Pi\cdot (\mathbf{x} \times \mathbf{x}') -
\Gamma\cdot (\mathbf{x} \times \mathbf{y}' - \mathbf{x}' \times
\mathbf{y})$ for any $\mathbf{x}, \mathbf{x}', \mathbf{y}, \mathbf{y}'
\in \mathbb{R}^3$. The heavy top equations $\dot{\Pi} = \Pi\times \Omega
+ Mg \ell \Gamma\times \chi$, $\dot{\Gamma} = \Gamma \times \Omega$
are Lie-Poisson equations on $\mathfrak{se}(3)^\ast$ for the Hamiltonian 
$h(\Pi, \Gamma) =  \frac{1}{2} \Pi \cdot \Omega + M g \ell \Gamma \cdot
\chi$ and the Lie-Poisson bracket $\{f, g\}(\Pi, \Gamma) = -\Pi\cdot
(\nabla_\Pi f \times \nabla_\Pi g)  - \Gamma\cdot (\nabla_\Pi f \times
\nabla_\Gamma g - \nabla_\Pi g \times \nabla_\Gamma f)$,
where $\nabla_\Pi$ and $\nabla_\Gamma $ denote the partial gradients.

Let $(\Pi(t), \Gamma(t)) $ be a periodic orbit of period $T $ of
the heavy top equations. After time $T$ by how much has the heavy top rotated in
space? The answer is provided by \eqref{general phase formula}:
\begin{align*}
\varphi &=  \frac{1}{\mu} \int\!\!\!\int_\mathcal{D} \omega_\mu + \frac{1}{\mu}\left(2h_\mu T
- 2 Mg \ell \int_0 ^T \Gamma(s) \cdot \chi  ds  \right)\\
&= \int\!\!\!\int_D \frac{2 \|\mathbb {I}\Gamma(s)\|^2 - (\Gamma(s)\cdot
\mathbb {I}\Gamma(s))(\operatorname{tr} \mathbb {I})}{(\Gamma(s)\cdot \mathbb
{I}\Gamma(s))^2} ds + \int_0^T \frac{ds}{\Gamma(s)\cdot \mathbb {I}\Gamma(s)} \,,
\end{align*}
where $D$ is the spherical cap on the unit sphere whose boundary is the closed curve
$\Gamma(t) $ and $\mathcal{D} $ is a two-dimensional submanifold of the orbit
$\mathcal{O}_ \mu $ bounded by the closed integral curve $(\Pi(t),
\Gamma(t))$. The first terms in each summand represent the geometric
phase and the second terms the dynamic phase.
\medskip

\noindent \textbf{Gauged Poisson structures.}
If the Lie group $G$ acts freely and properly on a smooth manifold
$Q$, then $(T^\ast Q)/G$ is a quotient Poisson manifold (see \cite{or1}),
where the quotient is taken relative  to the (left) lifted cotangent
action. The leaves of this Poisson manifold are the orbit reduced spaces
$\mathbf{J}^{-1}( \mathcal{O}_ \mu) /G$, where $\mathcal{O}_ \mu \subset
\mathfrak{g}^\ast$ is the coadjoint $G $-orbit through $\mu\in
\mathfrak{g}^\ast$ (see
\cite{or}). Is there an \textit{explicit\/} formula for this reduced
Poisson bracket on a manifold diffeomorphic to $(T^\ast Q)/G$? It turns
out that this question has two possible answers, once a connection on the
principal bundle $\pi: Q \rightarrow  Q/G $ is introduced. The discussion
below will also link to the fibration version of cotangent bundle
reduction.

In order to present these answers we review two bundle constructions. 
Let $G$ act freely and properly on the manifold $P$ and consider the  a
(left) principal $G$-bundle $\rho:P\rightarrow P/G: = M$. Let
$\tau: N\rightarrow M$ be a surjective submersion. Then the {\bfi
pull-back bundle\/} $\tilde{\rho} : (n, p) \in \tilde{P}: = \{(n,p) \in N
\times P \mid \rho(p)=\tau(n)\} \mapsto n \in N$ over $N$ is also a
principal (left) $G $-bundle relative to the action $g \cdot (n, p) : = 
(n, g \cdot p) $.

If there is a (left) $G$-action a manifold $V$, then the
diagonal $G $-action $g \cdot (p, v) = (g\cdot p, g\cdot v)$ on $P \times
V $ is also free and proper and one can form the {\bfi associated
bundle\/} $P \times_G V: =(P \times V)/G$ which is a locally trivial fiber
bundle $\rho_E: [p,v] \in E:= P \times_G V \mapsto \rho(p) \in M$ over
$M$ with fibers diffeomorphic to $V$. Analogously one can form the
associated fiber bundle  $\rho_{\tilde{E}}: \tilde{E} : = \tilde{P}
\times_G V \rightarrow N$.  Summarizing, the associated bundle
$\tilde{E}=\tilde{P}
\times_G V \rightarrow N$ is obtained from the principal bundle $\rho : P
\rightarrow M$, the surjective submersion $\tau : N \rightarrow M$, and
the $G$-manifold $V$ by pull-back and association, in this order.

These operations can be reversed. First form the associated bundle
$\rho_E: E = P \times_G V \rightarrow M$ and then 
pull it back by the surjective submersion $\tau:N \rightarrow M$ to $N $
to get the pull-back bundle $\tilde{\rho}_E: \tilde{E} \rightarrow N$. 
The map $\Phi:\tilde{P}\times_{G}V\rightarrow \tilde{E} $ defined by
$\Phi([(n,p),v]):=(n,[p,v])$ is an isomorphism of locally trivial fiber
bundles.

These general considerations will be used now to realize the quotient
Poisson manifold $(T^\ast Q)/G $ in two different ways. Let $Q$ be a
manifold and $G$  a Lie group (with Lie algebra
$\mathfrak{g}$) acting freely and properly on it. Let $A \in
\Omega^1(Q; \mathfrak{g})$ be a connection one-form on the left
$G$-principal bundle $ \pi : Q\rightarrow Q/G$. Pull back the $G $-bundle
$\pi: Q \rightarrow Q/G $ by the cotangent bundle projection
$\pi_{Q/G} :T^{\ast}(Q/G)\rightarrow Q/G$ to $T ^\ast(Q/G) $ to obtain the
$G$-principal bundle $\tilde{\pi}_{Q/G} : (\alpha_{[q]},q) \in \tilde{Q} :
=\{(\alpha_{[q]},q) \mid [q]=\pi(q), q \in Q\} \mapsto \alpha_{[q]}
\in T ^\ast(Q/G)$. This bundle is isomorphic to the annihilator $(VQ)^\circ
\subset T^\ast Q$  of the vertical bundle $VQ:= \ker T\pi \subset TQ$. 
Next, form the {\bfi coadjoint bundle\/}
$ \rho_S: S:=\tilde{Q}\times_{G}\mathfrak{g}^{\ast} \rightarrow
T^\ast(Q/G)$ of $\tilde{Q}$, $\rho_S((\alpha_{[q]},q), \mu)= 
\alpha_{[q]}$, that is, the associated vector bundle to the
$G$-principal bundle $\tilde{Q}\rightarrow T^\ast (Q/G)$ given by the
coadjoint representation of $G$ on $\mathfrak{g}^\ast$. The 
connection-dependent map $\Phi_{A}:S\rightarrow (T^{\ast}Q)/G $ defined
by $\Phi_{A}([(\alpha_{[q]},q),\mu]):
=[T^{\ast}_{q}\pi(\alpha_{[q]}) + A(q)^{\ast}\mu]$,
where $q \in Q $, $\alpha_q \in T ^\ast_q Q $, and $\mu\in \mathfrak{g}^\ast$, is
a vector bundle isomorphism over $Q/G$. The {\bfi Sternberg space\/}
is the Poisson manifold $(S, \{\cdot , \cdot \}_S)$, where $\{\cdot ,
\cdot \}_S$ is the pull-back to $S $ by $\Phi_A $ of the quotient Poisson
bracket on $(T ^\ast Q)/G$.

Next, we proceed in the opposite order. Construct first the coadjoint
bundle
$\rho_{\widetilde{\mathfrak{g}}^{\ast}} :[q, \mu] \in 
\widetilde{\mathfrak{g}}^{\ast} :=Q\times_{G}\mathfrak{g}^{\ast}
\mapsto [q] \in Q/G$ associated to the principal bundle
$\pi: Q \rightarrow Q/G$ and then pull it back by the cotangent bundle
projection $\pi_{Q/G}:T^{\ast}(Q/G)\rightarrow Q/G$ to $T ^\ast(Q/G)$ to
obtain the vector bundle $\rho_W: W: = \{(\alpha_{[q]},[q,\mu]) \mid
\pi_{Q/G}(\alpha_{[q]})=\rho_{\widetilde{\mathfrak{g}}^{\ast}}([q,\mu]) =
[q]\}$, $\rho_W(\alpha_{[q]},[q,\mu]) = \alpha_{[q]}$ over
$T^{\ast}(Q/G)$. Note that $W=T^{\ast}(Q/G)\oplus
{\widetilde{\mathfrak{g}}^{\ast}}$ and hence $W $ is also a vector bundle
over $Q/G $. Let $HQ$ be the horizontal subbundle defined by the connection
${A}$; thus
$TQ =HQ \oplus VQ$, where $H_q Q: = \ker{A}(q)$. For each $q \in
Q$ the linear map $T_q \pi |_{H_qQ}: H_qQ \rightarrow T_{[q]}(Q/G)$
is an isomorphism. Let $\operatorname{hor}_q\!:=
(T_q\pi|_{H_qQ})^{-1}\!:  T_{[q]}(Q/G)
\rightarrow H_q Q \subset T_q Q$ be the horizontal lift operator induced by
the connection ${A}$. Thus $\operatorname{hor}_{q}^{\ast} : T^\ast _q Q
\rightarrow T_{[q]}^\ast (Q/G)$ is a linear surjective map whose kernel is the
annihilator $(H_qQ)^\circ$ of the horizontal space.
The connection-dependent map $\Psi_{A}:(T ^\ast Q)/G \rightarrow W $
defined by $\Psi_{\mathcal{A}}([\alpha_{q}]):=
(\operatorname{hor}_{q}^{\ast}(\alpha_{q}), [q,\mathbf{J}(\alpha_{q})])$,
where $q \in Q $, $\alpha_q \in T ^\ast_q Q $, and $\mathbf{J}: T^\ast Q \rightarrow 
\mathfrak{g}^\ast$ is the momentum map of the lifted action, $\langle \mathbf{J}(\alpha_q), \xi
\rangle = \alpha_q((\xi_Q(q))$ for $\xi \in \mathfrak{g}$, is a vector
bundle isomorphism over $Q/G$ and $\Psi_A \circ \Phi_A = \Phi$.
The {\bfi Weinstein space\/} is the Poisson manifold $(W, \{\cdot , \cdot
\}_W)$, where $\{\cdot , \cdot\}_W $ is the push-forward by
$\Psi_{\mathcal{A}}$ of the Poisson bracket of $(T^\ast Q)/G$. In
particular, $\Phi:S \rightarrow W$ is a connection independent Poisson
diffeomorphism. The Poisson brackets  on $S $ and on $W $ are called
{\bfi gauged Poisson brackets}. They are expressed explicitly in terms of
various covariant derivatives induced on $S $ and on $W $ by the
connection $A \in \Omega ^1(Q; \mathfrak{g}) $.

Recall that the connection $A $ on the principal bundle $\pi: Q
\rightarrow Q/G$ naturally induces connections on pull-back bundles and
affine connections on  associated vector bundles. Thus both $S $ and $W $
carry covariant derivatives induced by $A $. They are given, according to
general definitions, in the cases under consideration, by:
\begin{itemize}
\item If $f \in C ^{\infty}(S) $, $s=[(\alpha_{[q]},q),\mu] \in S$, and
$v_{\alpha_{[q]}} \in T_{\alpha_{[q]}}T^\ast(Q/G)$, then
$\mathbf{d}_{\tilde{A}}^S  f(s)\in
T^{\ast}_{\alpha_{[q]}}T^{\ast}(Q/G)$ is defined  by
$\mathbf{d}_{\tilde{A}}^S f(s)\left(v_{\alpha_{[q]}}\right):
=\mathbf{d}f(s)
\left(T_{((\alpha_{[q]},q),\mu)}\pi_{\tilde{Q}\times\mathfrak{g}^{\ast}}\left(\left(v_{\alpha_{[q]}},
\mathrm{hor}_{q}(T_{\alpha_{[q]}}\tau (v_{\alpha_{[q]}})) \right),0 \right)
\right)$, where $\pi_{\tilde{Q}\times\mathfrak{g}^{\ast}} :
\tilde{Q}\times\mathfrak{g}^{\ast} \rightarrow
\tilde{Q}\times_G\mathfrak{g}^{\ast} = S$ is the orbit map. The symbol
$\mathbf{d}_{\tilde{A}}^S$ signifies that this is a covariant 
derivative on the associated bundle $S $ induced by the connection
$\tilde{A} $ on the principal $G $-pull-back bundle  $ \tilde{Q}
\rightarrow T ^\ast(Q/G) $. This connection $\tilde{A} $ is the pull-back
connection defined by $A $.
\item If $f \in C ^{\infty}(W)$,  $w=(\alpha_{[q]},[q,\mu])\in W$, and
$v_{\alpha_{[q]}} \in T_{\alpha_{[q]}}T^\ast (Q/G)$, then 
$\widetilde{\nabla}_{A}^W f(w) \in T^\ast _{\alpha_{[q]}} T^\ast
(Q/G) $ is defined by $\widetilde{\nabla}_{\mathcal{A}}^W f(w)
\left(v_{\alpha_{[q]}}\right) = \mathbf{d} f (w)\left(v_{\alpha_{[q]}}, 
T_{(q,\mu)}\pi_{Q\times \mathfrak{g}^\ast} 
\left( \operatorname{hor}_q (T_{\alpha_{[q]}}\tau_{Q/G} (v_{\alpha_{[q]}})),
0\right) \right)$, where $\pi_{Q \times \mathfrak{g}^\ast}: Q \times
\mathfrak{g}^\ast \rightarrow Q \times_G \mathfrak{g}^\ast =
\widetilde{\mathfrak{g}}^{\ast}$ is the orbit map. The symbol
$\widetilde{\nabla}_{\mathcal{A}}^W$ signifies that this is a covariant
derivative on the pull-back bundle $W $ induced by the covariant
derivative $\nabla_A $ on the coadjoint bundle
$\widetilde{\mathfrak{g}}^{\ast}$. This covariant  derivative $\nabla_A $
is induced on $\widetilde{\mathfrak{g}}^{\ast}$ by the connection $A $.
\item For $f\in C^{\infty}(W)$, we have
$\mathbf{d}_{\tilde{A}}^S(f\circ\Phi)
=(\widetilde{\nabla}_{A}^W f)\circ\Phi$.
\end{itemize}
To write the two gauged Poisson brackets on $S $ and on $W $ explicitly,
we denote by $\tilde{ \mathfrak{g}} = Q \times _G
\mathfrak{g}$ the adjoint  bundle of $\pi: Q \rightarrow Q/G $,
by $\Omega_{Q/G}
$ the  canonical symplectic  structure on
$T ^\ast(Q/G)$, by $B \in \Omega ^2(Q; \mathfrak{g}) $ the curvature of
$A$, and by $\mathcal{B}$ the $\tilde{\mathfrak{g}}$-valued two-form 
$\mathcal{B} \in \Omega ^2(Q/G; \tilde{\mathfrak{g}})$ on the base $Q/G$
defined by $\mathcal{B}([q])(u_{[q]},v_{[q]})
=[q,B(q)(u_q,v_q)]$, for any $u_q, v_q \in T_q Q$ that satisfy
$T_{q}\pi(u_{q})=u_{[q]}$ and $T_{q}\pi( v_{q})=v_{[q]}$. Note that both
$S ^\ast$ and $W ^\ast$ are Lie algebra bundles, that is, their fibers
are Lie algebras and the fiberwise Lie bracket operation depends smoothly
on the base point. If $f\in C^{\infty}(S)$, denote by $\delta f / \delta
s \in S ^\ast = \tilde{Q} \times _G \mathfrak{g}$ the usual fiber
derivative of $f $. Similarly, if $f \in C ^{\infty}(W) $ denote by
$\delta f/ \delta w \in W ^\ast$ the usual fiber derivative of $f $.
Finally, $\sharp:T^{\ast}(T^{\ast}(Q/G))\rightarrow T(T^{\ast}(Q/G))$ is
the vector  bundle isomorphism induced by $\Omega_{Q/G}$. The Poisson
bracket of $f, g \in C ^{\infty}(S) $ is given by
\begin{align*}
\{f,g\}_S(s) &= \Omega_{Q/G}(\alpha_{[q]})
\left(\mathbf{d}_{\tilde{\mathcal{A}}}^S
f(s)^{\sharp},\mathbf{d}_{\tilde{\mathcal{A}}}^S g(s)^{\sharp}\right)  
  -
\left\langle s,\left[\frac{\delta f}{\delta s},\frac{\delta g}{\delta s}\right]
\right\rangle \\
& \qquad +\left\langle v,
(\pi_{Q/G}^{\ast}\mathcal{B})(\alpha_{[q]})\left(\mathbf{d}_{\tilde{\mathcal{A}}}^S
f(s)^{\sharp},\mathbf{d}_{\tilde{\mathcal{A}}}^S
g(s)^{\sharp}\right)\right\rangle,
\end{align*}
where $v=[q,\mu] \in  \widetilde{\mathfrak{g}} ^\ast$. The Poisson
bracket $f, g \in C ^{\infty}(W) $ is given by
\begin{align*}
\{f,g\}_{W}(w) &= \Omega_{Q/G}(\alpha_{[q]})
\left(\widetilde{\nabla}_{\mathcal{A}}^W
f(w)^{\sharp},\widetilde{\nabla}_{\mathcal{A}}^W g(w)^{\sharp}\right)  
-\left\langle w,\left[\frac{\delta f}{\delta w},\frac{\delta g}{\delta
w}\right]\right\rangle \\
& \qquad +\left\langle v,
(\pi_{Q/G}^{\ast}\mathcal{B})(\alpha_{[q]})\left(\widetilde{\nabla}_{\mathcal{A}}^W
f(w)^{\sharp},
\widetilde{\nabla}_{\mathcal{A}}^W g(w)^{\sharp}\right)\right\rangle.
\end{align*}
Note that their structure is of the form: 
``canonical" bracket plus a (left) ``Lie-Poisson" bracket plus a
curvature coupling term.
\medskip

\noindent \textbf{The symplectic leaves of the Sternberg and
Weinstein spaces.} The map the map $\varphi_{A}:\tilde{Q} \times
\mathfrak{g}^\ast \rightarrow T^\ast Q$ given by
$\varphi_{A}\left(\left(\alpha_{[q]}, q \right), \mu\right):
=   T^{\ast}_{q}\pi(\alpha_{[q]}) +A(q)^{\ast}\mu$, where 
$\left(\left(\alpha_{[q]}, q \right), \mu\right) \in \tilde{Q} \times
\mathfrak{g}^\ast$, is a $G$-equivariant diffeomorphism; the $G$-action
on $T ^\ast Q $ is by cotangent lift and on  $\tilde{Q} \times
\mathfrak{g}^\ast $ is $g \cdot \left(\left(\alpha_{[q]}, q \right),
\mu\right) = \left(\left(\alpha_{[q]}, g \cdot q \right),
\operatorname{Ad}^\ast_{g ^{-1}}\mu\right)$. The pull back $\mathbf{J}_A$
of the momentum map to $\tilde{Q} \times \mathfrak{g}^\ast$ has the
expression $\mathbf{J}_A \left(\left(\alpha_{[q]}, q \right), \mu\right)
= \mu$, so if $\mathcal{O} \subset \mathfrak{g}^\ast$ is a coadjoint
orbit we have $\mathbf{J}_A ^{-1}(\mathcal{O})  = \tilde{Q} \times
\mathcal{O}$, and hence the orbit reduced manifold $\mathbf{J}_A ^{-1}(
\mathcal{O})/G $, whose connected components are the symplectic leaves of
$S$,  equals $\tilde{Q}\times_G\mathcal{O}$. Its symplectic form is the
{\bfi Sternberg minimal coupling form\/} 
$\tilde{\omega}_\mathcal{O}^- + \rho_S^\ast \Omega_{Q/G}$.

In this formula the two-form $\tilde{\omega}_\mathcal{O}^-$ has not been
defined yet. It is uniquely defined by the identity $\pi_{\tilde{Q}
\times \mathfrak{g}^\ast} \tilde{\omega}_\mathcal{O}^- = \mathbf{d}
\widehat{A} + \Pi_ \mathcal{O}\omega^- _ \mathcal{O}$, where $\omega^- _
\mathcal{O} $ is the minus orbit symplectic  form on $\mathcal{O}$ (see
\cite{or}), $\Pi_ \mathcal{O}: \tilde{Q} \times \mathcal{O} \rightarrow
\mathcal{O}$ is the projection on the second factor, and $\widehat{A}
\in \Omega ^2(\tilde{Q} \times \mathcal{O}) $ is the two-form given by
$\widehat{A}\left(( \alpha_{[q]}, q), \mu\right)
\left((u_{\alpha_{[q]}}, v_q), \nu \right) 
= - \left\langle \mu, A(q)(v_q) \right\rangle $ for $\left(( \alpha_{[q]}, q), \mu\right) \in \tilde{Q} \times \mathcal{O}$, 
$\left( u_{\alpha_{[q]}}, v_q\right) \in
T_{\left(\alpha_{[q]}, q \right)} \tilde{Q} $, and $\nu\in
\mathfrak{g}^\ast$. 

The symplectic leaves of the Weinstein  space $W $ are obtained by
pushing  forward by $\Phi$ the symplectic leaves of the Sternberg space.
They are the connected components of the symplectic manifolds
$\left(T^\ast(Q/G) \oplus (Q \times _G \mathcal{O}),
\Pi_{T^\ast(Q/G)}^\ast \Omega_{Q/G} + \Pi_{Q \times _G \mathcal{O}} ^\ast
\omega^-_{Q \times _G \mathcal{O}}\right)$, where
$\mathcal{O}$ is  a coadjoint orbit in  $\mathfrak{g}^\ast$, $\Omega_{Q/G} $
is the canonical symplectic form on $T ^\ast(Q/G) $, $\omega^-_{Q \times _G
\mathcal{O}}$ is a closed two-form on $Q \times _G \mathcal{O}$ to be
defined below, and
$\Pi_{T^\ast(Q/G)}: T^\ast(Q/G) \oplus (Q \times _G \mathcal{O}) \rightarrow 
T^\ast(Q/G)$, $\Pi_{Q \times _G \mathcal{O}}: T^\ast(Q/G) \oplus (Q \times _G
\mathcal{O}) \rightarrow Q \times _G \mathcal{O}$ are the projections.
The closed two-form $\omega^-_{Q \times _G \mathcal{O}} \in \Omega ^2
\left(Q \times _G \mathcal{O}\right)$ is uniquely determined by the
identity $\pi_{Q \times \mathcal{O}}^\ast \omega^-_{Q \times_G
\mathcal{O}} = \omega^-_{Q \times \mathcal{O}}$, where $\pi_{Q \times
\mathcal{O}} : Q \times \mathcal{O} \rightarrow Q \times _G \mathcal{O}$
is the orbit space projection, $\omega^-_{Q \times \mathcal{O}} \in
\Omega^2(Q \times \mathcal{O})$ is closed and given by $\omega^-_{Q \times \mathcal{O}}(q, \mu) \left((u_q, - \operatorname{ad}^\ast_
\xi\mu), (v_q,-\operatorname{ad}^\ast_ \eta \mu)\right) := 
-\mathbf{d} (\mathcal{A} \times \operatorname{id}_ \mathcal{O}) (q, \mu) 
\left((u_q, - \operatorname{ad}^\ast_\xi\mu), 
(v_q, -\operatorname{ad}^\ast_\eta \mu)\right) + \omega^-_ \mathcal{O}(
\mu) \left(\operatorname{ad}^\ast_\xi \mu, \operatorname{ad}^\ast_\eta \mu
\right)$, and $\mathcal{A} \times \operatorname{id}_ \mathcal{O} \in \Omega^1(Q \times
\mathfrak{g}^\ast)$ is given by $\left(\mathcal{A} \times
\operatorname{id}_\mathcal{O} \right) (q, \mu) \left(u_q,
-\operatorname{ad}^\ast_\xi\mu\right) = 
\left\langle \mu, \mathcal{A}(q)(u_q) \right\rangle$, for $q \in Q,
\mu\in \mathfrak{g}^\ast, u_q, v_q \in T_q Q, \xi,\eta\in \mathfrak{g}$. 

Thus on the Sternberg and Weinstein spaces both the Poisson bracket as
well as the symplectic form on the leaves have explicit connection
dependent formulas.

\end{document}